\documentclass{article}
\usepackage{latexsym}
\usepackage[varumlaut]{yfonts}
\usepackage{epsfig}
\usepackage{amsmath}
\usepackage{amssymb}
\usepackage{amsthm}
\usepackage{multirow}
\newtheorem{theorem}{Theorem}[section]
\newtheorem{corollary}[theorem]{Corollary}

\newtheorem{lemma}[theorem]{Lemma}

\theoremstyle{remark}

\theoremstyle{definition}
\newtheorem{definition}[theorem]{Definition}

\DeclareMathOperator\sign{sgn}

\DeclareMathOperator\Supp{supp}

\DeclareMathOperator\relint{relint}

\newcommand{\copos}[1]{\mathcal{C}_{#1}}
\newcommand{\NNM}[1]{\mathcal{N}_{#1}}
\newcommand{\SOZ}[1]{\mathcal{V}^{#1}}
\newcommand{\MSOZ}[1]{\mathcal{V}^{#1}_{\min}}

\begin{document}

\title{Minimal zeros of copositive matrices}

\author{Roland Hildebrand \thanks{%
WIAS, Mohrenstrasse 39, 10117 Berlin, Germany
({\tt roland.hildebrand@wias-berlin.de}).}}

\maketitle

\begin{abstract}
Let $A$ be an element of the copositive cone $\copos{n}$. A zero $u$ of $A$ is a nonzero nonnegative vector such that $u^TAu = 0$. The support of $u$ is the index set $\Supp{u} \subset \{1,\dots,n\}$ corresponding to the positive entries of $u$. A zero $u$ of $A$ is called minimal if there does not exist another zero $v$ of $A$ such that its support $\Supp{v}$ is a strict subset of $\Supp{u}$. We investigate the properties of minimal zeros of copositive matrices and their supports. Special attention is devoted to copositive matrices which are irreducible with respect to the cone $S_+(n)$ of positive semi-definite matrices, i.e., matrices which cannot be written as a sum of a copositive and a nonzero positive semi-definite matrix. We give a necessary and sufficient condition for irreducibility of a matrix $A$ with respect to $S_+(n)$ in terms of its minimal zeros. A similar condition is given for the irreducibility with respect to the cone $\NNM{n}$ of entry-wise nonnegative matrices. For $n = 5$ matrices which are irreducible with respect to both $S_+(5)$ and $\NNM{5}$ are extremal. For $n = 6$ a list of candidate combinations of supports of minimal zeros which an exceptional extremal matrix can have is provided.
\end{abstract}



\section{Introduction}

A real symmetric $n \times n$ matrix $A$ is called {\it copositive} if $x^TAx \geq 0$ for all $x \in \mathbb R_+^n$. The set of copositive matrices forms a convex cone, the {\it copositive cone} $\copos{n}$. This matrix cone is of interest for combinatorial optimization, for surveys see~\cite{Duer10,HUS10,BSU12}. However, verifying copositivity of a given matrix is a co-NP-complete problem~\cite{MurtyKabadi87}. It is a classical result by Diananda \cite[Theorem 2]{Diananda62} that for $n \leq 4$ the copositive cone can be described as the sum of the cone of positive semi-definite matrices $S_+(n)$ and the cone of element-wise nonnegative symmetric matrices $\NNM{n}$. In general, this sum is a subset of the copositive cone, $S_+(n) + \NNM{n} \subset \copos{n}$. Horn showed that for $n \geq 5$ the inclusion is strict \cite[p.25]{Diananda62}.

A nonzero vector $u \in \mathbb R_+^n$ is called a {\it zero} of a copositive matrix $A$ if $u^TAu = 0$. It has been recognised early that the zero set of a copositive matrix is a useful tool in the study of the structure of the cone $\copos{n}$ \cite{Diananda62,HallNewman63}. In \cite{Baumert67} Baumert considered the possible zero sets of matrices in ${\cal C}_5$. He provided a partial classification of the zero sets of matrices $A \in \copos{5}$ which are {\it irreducible} with respect to the cone $\NNM{5}$, i.e., which cannot be written as a nontrivial sum $A = C + N$, where $C$ is copositive and $N$ is element-wise nonnegative. In \cite{DDGH13a} this classification was completed and a necessary and sufficient condition for irreducibility of a copositive matrix $A \in \copos{n}$ with respect to the cone $\NNM{n}$ was given in terms of its zero set. This allowed the classification of the extreme rays of the cone ${\cal C}_5$ in \cite{Hildebrand12a}.

In \cite{Baumert65} Baumert introduced the concept of maximal zeros of a copositive matrix. He called a zero $u$ of $A$ {\it maximal} if for no other zero $v$ of $A$, the index set of positive entries of $u$ is a strict subset of the index set of positive entries of $v$. In this note we introduce and investigate the concept of minimal zeros of copositive matrices. Here a zero $u$ of $A$ is called {\it minimal} if for no other zero $v$ of $A$, the index set of positive entries of $v$ is a strict subset of the index set of positive entries of $u$. We consider some properties of the set of minimal zeros of a copositive matrix and derive necessary and sufficient conditions for a copositive matrix $A$ to be irreducible with respect to the cone $\NNM{n}$ or $S_+(n)$ in terms of its set of minimal zeros. In contrast to maximal zeros, or zeros of copositive matrices in general, a minimal zero is determined up to scaling by a positive constant by the index set of its positive entries. Thus a copositive matrix can essentially have only a finite number of minimal zeros, which opens the way to a combinatorial approach.

The obtained results can potentially be used in order to obtain a classification of the extreme rays of the cone $\copos{n}$ for small $n$. In application to the case $n = 5$, we show that a matrix $A \in \copos{5}$ with positive diagonal elements is irreducible with respect to both cones $S_+(5)$ and $\NNM{5}$ if and only if its set of minimal zeros is one of two types. These types correspond to the two types of exceptional extreme rays of $\copos{5}$ which have been obtained in \cite{Hildebrand12a}, i.e., extreme rays which are not contained in the sum $S_+(5) + \NNM{5}$ \cite{JohnsonReams08}. Thus, by using the results of this paper, the classification of the extreme rays of $\copos{5}$ can be reduced to the consideration of just two cases, in contrast to the approximately 30 cases which have been considered in \cite{Baumert67,DDGH13a} and on which the classification in \cite{Hildebrand12a} is based. For $n = 6$, a matrix $A \in \copos{6}$ with positive diagonal elements is irreducible with respect to both $S_+(6)$ and $\NNM{6}$ if and only if its set of minimal zeros is one of 44 types.

The remainder of the paper is structured as follows. In the next section we provide necessary definitions and collect some results from the literature for later use. In Section \ref{sec:minimal}, we characterize minimal zeros in different ways and establish conditions on the combinations of minimal zeros that a copositive matrix can have. In Section \ref{sec:irred} we consider irreducibility of a copositive matrix with respect to the cones of positive semi-definite and nonnegative matrices, respectively. In Section \ref{sec:exceptional} we apply the results in order to restrict the combinations of minimal zeros that can occur in exceptional extreme copositive matrices. We provide a list of combinations for the cone $\copos{6}$. Finally, we give a summary in the last section.

\section{Notations and preliminaries}

We shall denote vectors with lower-case letters and matrices with upper-case letters. Individual entries of a vector $u$ or a matrix $A$ will be denoted by $u_i$, $A_{ij}$, respectively. For a matrix $A$ and a vector $u$ of compatible size, the $i$-th element of the vector $Au$ will be denoted by $(Au)_i$. Inequalities $u \geq 0$ on vectors will be meant element-wise. We denote by ${\bf 1} = (1,\dots,1)^T$ the all-ones vector. Let further $E_{ij}$ be the $n \times n$ matrix that has zero entries everywhere except at $(i,j)$ and $(j,i)$, where it has entries 1.

For a subset $I \subset \{1,\dots,n\}$ we denote by $A_I$ the principal submatrix of $A$ whose elements have row and column indices in $I$, i.e. $A_I = (A_{ij})_{i,j \in I}$. Similarly for a vector $u \in \mathbb R^n$ we define the subvector $u_I = (u_i)_{i \in I}$.

We call a nonzero vector $u \in \mathbb R_+^n$ a {\it zero} of a copositive matrix $A \in \copos{n}$ if $u^TAu = 0$. We denote the set of zeros of $A$ by $\SOZ{A} = \{u \in \mathbb R^n_+\setminus\{0\}\mid u^TAu = 0 \}$. For a vector $u \in \mathbb R^n$ we define its {\it support} as $\Supp{u} = \{i\in\{1,\ldots,n\}\mid u_i \not= 0 \}$. A zero $u$ of a copositive matrix $A$ is called {\it minimal} if there exists no zero $v$ of $A$ such that the inclusion $\Supp{v} \subset \Supp{u}$ holds strictly. We shall denote the set of minimal zeros of a copositive matrix $A$ by $\MSOZ{A}$. The {\it support set} of $A$ is the set $\Supp{\SOZ{A}} = \{ \Supp{u} \,|\, u \in \SOZ{A} \}$, and the {\it minimal support set} is the set $\Supp{\MSOZ{A}} = \{ \Supp{u} \,|\, u \in \MSOZ{A} \}$.

An element $A \in \copos{n}$ is called {\it extremal} if the conditions $A = B + C$, $B,C \in \copos{n}$ imply the existence of nonnegative numbers $\lambda,\mu$ such that $B = \lambda A$, $C = \mu A$. The conic hull of a nonzero extremal element $A \in \copos{n}$ is an {\it extreme ray}. Following \cite{JohnsonReams08}, if $A \not\in S_+(n) + \NNM{n}$, then $A$ and the extreme ray it generates are called {\it exceptional}. %We will, however, use the term {\it exceptional} for any copositive matrix $A \not\in S_+(n) + \NNM{n}$, not necessarily an extremal one.

\begin{definition} \cite[Definition 1.1]{DDGH13a} \label{def_irred}
For a matrix $A\in\copos{n}$ and a subset $\mathcal{M} \subset \copos{n}$, we say that $A$ is {\sl irreducible with respect to $\mathcal{M}$} if there do not exist $\gamma>0$ and $M\in\mathcal{M}\setminus\{0\}$ such that $A-\gamma M\in\copos{n}$.
\end{definition}

Note that this definition differs from the concept of an irreducible matrix that is normally used in matrix theory. For simplicity we speak about irreducibility with respect to $M$ when $\mathcal{M}=\{M\}$. In our paper, we shall be concerned with the cases
\[ \mathcal{M} = S_+(n),\quad\mathcal{M} = \{ww^T\},\quad\mathcal{M} = \NNM{n},\quad \mbox{and}\quad \mathcal{M} = \{E_{ij}\}.
\]

{\lemma \label{basic_irred} Let $A \in \copos{n}$ and $\mathcal{M} \subset \copos{n}$. Then the following are equivalent.

(a) $A$ is irreducible with respect to $\mathcal{M}$,

(b) $A$ is irreducible with respect to $M$ for all $M \in \mathcal{M}$,

(c) $A$ is irreducible with respect to $\mathbb R_+\mathcal{M}$,

(d) $A$ is irreducible with respect to the convex conic hull of $\mathcal{M}$. }

\begin{proof}
The equivalence of (a)--(c) and the implication (d) $\Rightarrow$ (a) follow directly from Definition \ref{def_irred}. Let us show the implication (a) $\Rightarrow$ (d).

For the sake of contradiction, assume (a) and let $M = \sum_{k=1}^m \alpha_k M_k$ be a nonzero element of the convex conic hull of $\mathcal{M}$, with $\alpha_k > 0$ and $M_k \in \mathcal{M}\setminus\{0\}$ for all $k$, such that $A - \gamma M \in \copos{n}$ for some $\gamma > 0$. We then also have $(A - \gamma M) + \gamma\sum_{k=2}^m \alpha_k M_k = A - \gamma\alpha_1 M_1 \in \copos{n}$, as this is a sum of copositive matrices. But this contradicts (a), because $\gamma\alpha_1 > 0$ and $M_1 \in \mathcal{M}\setminus\{0\}$.
\end{proof}

In particular, $A \in \copos{n}$ is irreducible with respect to $S_+(n)$ if and only if it is irreducible with respect to $ww^T$ for every nonzero vector $w \in \mathbb R^n$, and it is irreducible with respect to $\NNM{n}$ if and only if it is irreducible with respect to $E_{ij}$ for all $i,j = 1,\dots,n$.

Note that if a matrix $A$ is on an exceptional extreme ray of $\copos{n}$, then $A$ must be irreducible with respect to both $S_+(n)$ and $\NNM{n}$.

\medskip

Finally we collect some results from the literature that will be used later on.

{\lemma \cite[Lemma 2.4]{DDGH13a} \label{zeroPSD} Let $A \in \copos{n}$ and $u \in \SOZ{A}$. Then the principal submatrix $A_{\Supp{u}}$ is positive semi-definite. }

{\lemma \cite[Lemma 2.5]{DDGH13a} \label{zero2PSD} Let $A \in \copos{n}$ and $u \in \SOZ{A}$. Then $(Au)_i = 0$ for all $i \in \Supp(u)$. }

{\lemma \cite[Theorem 7.2]{Dickinson11} \label{Th72D} Let $A \in \copos{n}$ be a copositive matrix and $I \subset \{1,\dots,n\}$ an index set of cardinality $k$. Suppose that the principal submatrix $A_I$ is positive semi-definite, and let $u \in \mathbb R^n$ be a nonzero vector such that $\Supp{u} \subset I$. Then $u \in \SOZ{A}$ if and only if $u_I \in \mathbb R_+^k \cap \ker A_I$. }

{\lemma \cite[p.200]{Baumert66} \label{first_order} Let $A \in \copos{n}$ and $u \in \SOZ{A}$. Then $Au \geq 0$. }

{\lemma \cite[Theorem 2.6]{DDGH13a} \label{irreducibility_old} Let $A \in \copos{n}$, and let $i,j \in \{1,\dots,n\}$. Then $A$ is irreducible with respect to $E_{ij}$ if and only if there exists a zero $u$ of $A$ such that $(Au)_i = (Au)_j = 0$ and $u_i + u_j > 0$. }

In \cite{DDGH13a} the lemma was stated for $n \geq 2$, but it is easily seen that the assertion holds also for $n = 1$.

{\lemma \cite[Corollary 4.4]{DDGH13a} \label{two_zeros} Let $A \in \copos{n}$ with $A_{ii} = 1$ for all $i$, and let $u \in \SOZ{A}$ with $|\Supp{u}| = 2$. Then the two positive elements of $u$ are equal. }

{\lemma \cite[Corollary 4.14]{DDGH13a} \label{supp_nminus2} Let $A \in \copos{n}$ be irreducible with respect to $\NNM{n}$. If there exists $u \in \SOZ{A}$ with $|\Supp{u}| \geq n-1$, then $A \in S_+(n)$. }

{\lemma \cite[Lemma 4.7]{DDGH13a} \label{3zero} Let $A \in \copos{3}$ with $A_{ii} = 1$ for $i = 1,2,3$. Then $A$ is irreducible with respect to $\{E_{12},E_{13},E_{23}\}$ if and only if it is of the form $A = \begin{pmatrix} 1 & -\cos\phi_3 & -\cos\phi_2 \\ -\cos\phi_3 & 1 & -\cos\phi_1 \\ -\cos\phi_2 & -\cos\phi_1 & 1 \end{pmatrix}$ for some scalars $\phi_1,\phi_2,\phi_3 \geq 0$ satisfying $\phi_1 + \phi_2 + \phi_3 = \pi$. }

{\lemma \cite[Theorem 5.6]{DDGH13a} \label{ineq4} Let $A \in \copos{5}$ with $A_{ii} = 1$ for all $i = 1,\dots,5$ be irreducible with respect to $E_{ij}$ for all $1 \leq i < j \leq 5$. Then either $A \in S_+(5)$, or there exists a permutation matrix $P$ and scalars $\theta_1,\dots,\theta_5 \geq 0$ such that $\sum_{i=1}^5 \theta_i < \pi$ and
\begin{equation} \label{Smatrix}
PAP^T = \begin{pmatrix}
       1 & -\cos\theta_1 & \cos(\theta_1+\theta_2) & \cos(\theta_4+\theta_5) & -\cos\theta_5 \\
       -\cos\theta_1 & 1 & -\cos\theta_2 & \cos(\theta_2+\theta_3) & \cos(\theta_5+\theta_1) \\
       \cos(\theta_1+\theta_2) & -\cos\theta_2 & 1 & -\cos\theta_3 & \cos(\theta_3+\theta_4) \\
       \cos(\theta_4+\theta_5) & \cos(\theta_2+\theta_3) & -\cos\theta_3 & 1 & -\cos\theta_4 \\
        -\cos\theta_5 & \cos(\theta_5+\theta_1) & \cos(\theta_3+\theta_4) & -\cos\theta_4 & 1
   \end{pmatrix}.
\end{equation} }

Finally we provide a result that is closely linked with the semi-definite approximation of the MAXCUT problem by Goemans and Williamson \cite{GoemansWilliamson}.

\begin{definition} \label{def_maxcut}
The {\sl MAXCUT polytope} ${\cal MC}_n \subset S_+(n)$ is the convex hull of all matrices $A \in S_+(n)$ such that $A_{ij} \in \{-1,+1\}$ for all $i,j = 1,\dots,n$, i.e., all matrices of the form $vv^T$, $v \in \{-1,+1\}^n$.
\end{definition}

The following lemma is a consequence of \cite[Lemma 3.2]{GoemansWilliamson}.

{\lemma \cite[Corollary 4.3]{Hirschfeld} \label{lem_maxcut} Let $A \in S_+(n)$ be a positive semi-definite matrix with $A_{ii} = 1$, $i = 1,\dots,n$. Let $B$ be the real symmetric $n \times n$ matrix defined entry-wise by $B_{ij} = \frac{2}{\pi}\arcsin A_{ij}$, $i,j = 1,\dots,n$. Then $B \in {\cal MC}_n$. }

\section{Minimal zeros of copositive matrices} \label{sec:minimal}

In this section we consider properties of minimal zeros of general copositive matrices. First we state an auxiliary result.

{\lemma \label{zero_cone} Let $A$ be a copositive matrix and $u \in \SOZ{A}$ with support $\Supp{u} = I$. Let $k$ be the cardinality of $I$ and denote the intersection $\mathbb R_+^k \cap \ker A_I$ by $K$. Let $v \in \mathbb R^n$ be a nonzero vector such that $\Supp{v} \subset I$. Then the following are equivalent.

(a) $v \in \SOZ{A}$,

(b) $v_I \in K$. }

\begin{proof}
By Lemma \ref{zeroPSD} the principal submatrix $A_I$ is positive semi-definite. The proof of the lemma now follows from Lemma \ref{Th72D}.
%For $w \in \mathbb R^k$, we hence have $w^TA_Iw = 0$ if and only if $w \in \ker A_I$.
%
%Assume that (a) holds. Then by definition $v \not= 0$, $\Supp{v} \subset I$, $v_I \geq 0$, and $v^TAv = 0$. The last relation implies $v_I^TA_Iv_I = 0$, and by the above $v_I \in \ker A_I$. Hence $v_I \in K$, which proves (b).
%
%Let now (b) hold. Then $v \geq 0$ and $v^TAv = v_I^TA_Iv_I = 0$, which implies $v \in \SOZ{A}$. Moreover, every positive entry of $v$ has index in $I$, and $\Supp{v} \subset I$. This shows (a).
\end{proof}

The lemma states that the set of zeros of a copositive matrix $A$ whose support is contained in the support of some fixed zero, is a convex polyhedral cone. This does not hold for the set of all zeros, which is not convex in general. We now relate the minimal zeros to the extreme rays of this cone. We first characterize these extreme rays.

{\lemma \label{char_extreme} Let $L \subset \mathbb R^k$ be a nonempty linear subspace, let $K = \mathbb R_+^k \cap L$, and let $u \in K$ be a nonzero vector. Then the following are equivalent.

(a) $u$ is an extremal element of $K$,

(b) if $v \in K$ and $\Supp{v} \subset \Supp{u}$, then $v$ is a multiple of $u$. }

\begin{proof}
Let $u$ be extremal and let $v \in K$ be such that $\Supp{v} \subset \Supp{u}$. Then there exists $\varepsilon > 0$ such that $w = u - \varepsilon v \geq 0$. Since $u,v \in L$, we also have $w \in K$. Then $u = w + \varepsilon v$, and by extremality of $u$ the vectors $v,w$ must be multiples of $u$. This proves the implication (a) $\Rightarrow$ (b).

Let us now assume (b) and suppose that $u = v + w$ for some $v,w \in K$. Since $v,w \geq 0$, we have $\Supp{v},\Supp{w} \subset \Supp{u}$. By condition (b), $v,w$ are then multiples of $u$. This proves the extremality of $u$.
\end{proof}

{\lemma \label{minimal_extremal} Assume the conditions of Lemma \ref{zero_cone}. Then the following are equivalent.

(a) $v$ is a minimal zero of $A$,

(b) $v_I$ is an extremal element of $K$. }

\begin{proof}
For any vector $y \in \mathbb R^k$, let $\tilde y \in \mathbb R^n$ with $\Supp{\tilde y} \subset I$ be defined by $\tilde y_I = y$. By Lemma \ref{zero_cone}, for every nonzero vector $y \in K$ the vector $\tilde y$ is a zero of $A$.

Assume condition (a). Then condition (a) of Lemma \ref{zero_cone} holds, and hence also condition (b) of this lemma. Assume for the sake of contradiction that $v_I$ is not an extremal element of $K$. Then there exist linearly independent vectors $w,z \in K$ such that $v_I = \frac{w + z}{2}$. Consider the proper affine line $\{ y(\lambda) = \lambda w + (1-\lambda)z \,|\, \lambda \in \mathbb R \}$ in $\mathbb R^k$ and the corresponding proper affine line $\{ \tilde y(\lambda) = \lambda \tilde w + (1-\lambda)\tilde z \,|\, \lambda \in \mathbb R \}$ in $\mathbb R^n$. Define the interval $J = \{ \lambda \in \mathbb R \,|\, y(\lambda) \in K \}$. Then $\tilde y(\lambda) \in \SOZ{A}$ for all $\lambda \in J$. By closedness of $K$ this interval is closed, and by $w,z \in K$ we have $[0,1] \subset J$. Since $w,z \geq 0$ and $v_I = \frac{w + z}{2}$, we have $\Supp{w},\Supp{z} \subset \Supp{v_I}$. Hence the indices of the nonzero elements of $y(\lambda)$ are contained in $\Supp{v_I}$ for every $\lambda \in \mathbb R$. In particular, for every $\lambda \in J$ we have $\Supp{y(\lambda)} \subset \Supp{v_I}$ and $\Supp{\tilde y(\lambda)} \subset \Supp{v}$. By the minimality of $v$, we then have $\Supp{\tilde y(\lambda)} = \Supp{v}$ for every $\lambda \in J$, and hence also $\Supp{y(\lambda)} = \Supp{v_I}$. But since the indices of the nonzero elements of $w-z$ are contained in $\Supp{v_I}$, the set $\{ \lambda \in J \,|\, \Supp{y(\lambda)} = \Supp{v_I} \}$ must be open. It follows that $J = \mathbb R$, which contradicts the pointedness of the cone $K$. This proves (b).

Let us now assume (b). Then condition (b) of Lemma \ref{zero_cone} holds, and hence also condition (a) of this lemma. Assume for the sake of contradiction that there exists a zero $w$ of $A$ such that $\Supp{w} \subset \Supp{v}$ strictly. Then we have also $\Supp{w_I} \subset \Supp{v_I}$ and $w_I \in K$. By Lemma \ref{char_extreme} the extremality of $v_I$ implies that $w_I$ is a multiple of $v_I$. But then $w$ is a multiple of $v$, which contradicts the strictness of the inclusion $\Supp{w} \subset \Supp{v}$. This proves (a).
\end{proof}

{\corollary \label{zero_decomp} Let $A$ be a copositive matrix and $u \in \SOZ{A}$. Then $u$ can be represented as a finite sum of minimal zeros of $A$. }

\begin{proof}
Let $I = \Supp{u}$ and let the cone $K$ be defined as in Lemma \ref{zero_cone}. By this lemma, we have $u_I \in K$. Then there exists a finite number of nonzero extremal elements $\tilde v,\dots,\tilde w$ of $K$ such that $\tilde v + \dots + \tilde w = u_I$. Define vectors $v,\dots,w \in \mathbb R^n$ such that $\Supp{\tilde v},\dots,\Supp{\tilde w} \subset I$ and $v_I = \tilde v,\dots,w_I = \tilde w$. Then by Lemma \ref{minimal_extremal}, $v,\dots,w$ are minimal zeros of $A$, and by construction $v + \dots + w = u$.
\end{proof}

Next we show that up to multiplication by a constant, a minimal zero is defined by its support.

{\lemma \label{uniqueness} Let $A$ be a copositive matrix and $u \in \SOZ{A}$. Then the following are equivalent.

(a) $u$ is a minimal zero of $A$,

(b) if $v$ is another zero of $A$ with support $\Supp{v} \subset \Supp{u}$, then there exists $\mu > 0$ such that $v = \mu u$. }

\begin{proof}
Let the cone $K$ be defined as in Lemma \ref{zero_cone}.

Assume condition (a). Then by Lemma \ref{minimal_extremal} the vector $u_I$ is an extremal element of $K$. Note that for every vector $w \in K$ we have $\Supp{w} \subset \Supp{u_I}$, and hence by Lemma \ref{char_extreme} the cone $K$ is 1-dimensional. Let now $v \in \SOZ{A}$ with $\Supp{v} \subset \Supp{u}$. By Lemma \ref{zero_cone} we then have $v_I \in K$, and by the preceding $v_I$ is a multiple of $u_I$. It follows that $v$ is a multiple of $u$. Condition (b) now easily follows.

Assume condition (b). Then for every $v \in \SOZ{A}$ with $\Supp{v} \subset \Supp{u}$ we have $\Supp{v} = \Supp{u}$. Hence $u$ is a minimal zero by definition, which proves (a).
%Let $w \in K$ be a nonzero vector, and define $\tilde w \in \mathbb R^n$ such that $\Supp{\tilde w} \subset I$ and $\tilde w_I = w$. By Lemma \ref{zero_cone} we have $\tilde w \in \SOZ{A}$ with $\Supp{\tilde w} \subset \Supp{u}$. By condition (b) $\tilde w$ is a multiple of $u$, and hence $w$ is a multiple of $u_I$. It follows that $K$ is 1-dimensional, and $u_I$ is an extremal element. By Lemma \ref{minimal_extremal} $u$ is then a minimal zero, which proves (a).
\end{proof}

{\corollary Let $A$ be a copositive matrix. Then the number of equivalence classes of minimal zeros of $A$ with respect to multiplication by a positive constant is finite. \qed }

\medskip

The classes of minimal zeros are hence in a one-to-one correspondence with the minimal support set $\Supp{\MSOZ{A}}$. They are also in a one-to-one correspondence with the elements of the sets ${\cal X}_i$ from \cite[Method 7.3]{Dickinson11}, which relate to the extreme rays of the intersection of the nonnegative orthant with the kernel of the maximal positive semi-definite principal submatrices of $A$.

Next we give a characterization of minimal zeros in terms of principal submatrices.

{\lemma \label{corank1} Let $A \in \copos{n}$ be a copositive matrix and let $I \subset \{1,\dots,n\}$ be a nonempty index set. Then the following are equivalent.

(a) $A$ has a minimal zero with support $I$,

(b) the principal submatrix $A_I$ is positive semi-definite with corank 1, and the generator of the kernel of $A_I$ can be chosen such that all its elements are positive. }

\begin{proof}
Assume the notations of Lemma \ref{zero_cone}.

Assume condition (a), and let $u$ be the minimal zero. Then by Lemma \ref{zeroPSD} the submatrix $A_I$ is positive semi-definite. The vector $u_I$ is in the interior of $\mathbb R_+^k$ and is by Lemma \ref{minimal_extremal} an extremal element of the cone $K = \mathbb R_+^k \cap \ker A_I$. It follows that $K$ is 1-dimensional, and hence $\ker A_I$ is 1-dimensional and generated by $u_I$. This proves (b).

Assume condition (b). Choose a vector $u \in \mathbb R_+^n$ such that $\Supp{u} = I$ and $u_I$ generates the kernel of $A_I$. Then $u_I \in K$. Since the kernel of $A_I$ is 1-dimensional by assumption, the vector $u_I$ is also an extremal element of $K$. By Lemma \ref{minimal_extremal} it follows that $u$ is a minimal zero, which proves (a).
\end{proof}

{\corollary \label{cor:corank1} Let $A \in \copos{n}$ be a copositive matrix and $u$ a minimal zero of $A$ with support $I = \Supp{u}$. Then for every proper subset $J \subset I$ we have that the principal submatrix $A_J$ is positive definite. }

\begin{proof}
By Lemma \ref{corank1} the principal submatrix $A_I$ is positive semi-definite, has corank 1, and its kernel is generated by a vector all whose elements are nonzero. Therefore $A_J$ is both positive semi-definite and non-degenerate, which implies that it is positive definite.
\end{proof}

Next we shall consider $m$-tuples of minimal zeros with overlapping supports. We begin with an auxiliary result.

{\lemma \label{PDlemma} Let $A \in \copos{n}$ be a copositive matrix and $I \subset \{1,\dots,n\}$ an index set such that $A_I$ is positive definite. Let $u \in \SOZ{A}$ such that $(\Supp{u}) \setminus I = \{k\}$ consists of exactly one element, and let $u$ be normalized such that $u_k = 1$. Then we have $A_{kk} = u_I^TA_Iu_I$, and $u$ is a minimal zero. }

\begin{proof}
By Lemma \ref{zeroPSD} the principal submatrix $A_{\Supp{u}}$ is positive semi-definite. It has the kernel vector $u_{\Supp{u}}$, which consists of positive elements only, and its principal submatrix $A_{(\Supp{u}) \setminus \{k\}}$ is positive definite. Hence $A_{\Supp{u}}$ is of corank 1 and by Lemma \ref{corank1} $A$ has a minimal zero $v$ with support $\Supp{u}$. By Lemma \ref{uniqueness} $u$ is proportional to $v$ and hence minimal.

By Lemma \ref{zero2PSD} we have $A_{\Supp{u}}u_{\Supp{u}} = 0$. It follows that $(Au)_i = \sum_{j \in I} A_{ij}u_j + A_{ik} = 0$ for every $i \in \Supp{u}$. In particular, $A_{ik} = -\sum_{j \in I} A_{ij}u_j$ for all $i \in \Supp{u}$. Setting $i = k$ yields $A_{kk} = -\sum_{l \in I} A_{lk}u_l = \sum_{j,l \in I} A_{lj}u_ju_l = u_I^TA_Iu_I$.
\end{proof}

{\corollary \label{PDcor} Let $A \in \copos{n}$ be a copositive matrix and $I \subset \{1,\dots,n\}$ an index set such that $A_I$ is positive definite. Let $u,v$ be zeros of $A$ such that $(\Supp{u}) \setminus I = \{k\}$, $(\Supp{v}) \setminus I = \{l\}$ for some indices $k,l$, and let $u,v$ be normalized such that $u_k = v_l = 1$. Then $(Au)_l = A_{kl} + \sum_{j \in I} u_j(Av)_j - v_I^TA_Iu_I$. If $u \not= v$ after normalization, then $k \not= l$. }

\begin{proof}
We have $(Au)_l = \sum_{j \in I} A_{lj}u_j + A_{kl} = \sum_{j \in I} ((Av)_j - \sum_{l \in I}A_{jl}v_l)u_j + A_{kl}$, which proves our first claim.

Let now $k = l$. We then have
\begin{eqnarray*}
0 &=& u^TAu = u_I^TA_Iu_I + 2\sum_{i \in I}A_{il}u_i + A_{ll} \\ &=& u_I^TA_Iu_I + 2\sum_{i \in I}\left((Av)_i - \sum_{j \in I}A_{ij}v_j\right)u_i + v_I^TA_Iv_I \\
&=& u_I^TA_Iu_I + 2\sum_{i \in I}(Av)_iu_i - 2u_I^TA_Iv_I + v_I^TA_Iv_I \geq (u_I-v_I)^TA_I(u_I-v_I).
\end{eqnarray*}
Here we used Lemma \ref{PDlemma} for the second relation. Since $A_I \succ 0$, it follows that $u_I = v_I$ and hence $u = v$. The second claim of the corollary now easily follows.
\end{proof}

{\theorem \label{PDth} Let $A \in \copos{n}$ be a copositive matrix and $I \subset \{1,\dots,n\}$ an index set such that the principal submatrix $A_I$ is positive definite. Let $u^1,\dots,u^m$ be zeros of $A$ such that $(\Supp{u^l}) \setminus I = \{k^l\}$ consists of exactly one element, and let $u^l$ be normalized such that $u^l_{k^l} = 1$, $l = 1,\dots,m$. Suppose that the zeros $u^1,\dots,u^m$ are mutually different after normalization. Suppose further that $\Supp{u^r_I} \subset \Supp{u^{r+1}_I}$ for all $r = 1,\dots,m-1$.

Then the indices $k^1,\dots,k^m$ are mutually different, and $u^1,\dots,u^m$ are minimal zeros. Moreover, if $v \in \SOZ{A}$ is a zero satisfying $\Supp{v} \subset I \cup \{k^1,\dots,k^m\}$, then $v = \sum_{i=1}^m \alpha_iu^i$ for some nonnegative scalars $\alpha^i$. If in addition $v$ is minimal, then there exists $l \in \{1,\dots,m\}$ and $\alpha > 0$ such that $v = \alpha u^l$. }

\begin{proof}
By Lemma \ref{PDlemma} the zeros $u^1,\dots,u^m$ are minimal, and by Corollary \ref{PDcor} the indices $k^1,\dots,k^m$ are mutually different.

For the sake of notational simplicity, let us assume without loss of generality that $k^r = r$ for $r = 1,\dots,m$.

Consider $r,s \in \{1,\dots,m\}$ such that $r \geq s$. By Lemma \ref{zero2PSD} we have $(Au^r)_i = 0$ for all $i \in I$ such that $u^r_i > 0$. But then also $(Au^r)_i = 0$ for all $i \in I$ such that $u^s_i > 0$, because $\Supp{u^s_I} \subset \Supp{u^r_I}$. Hence $(Au^r)_iu^s_i = 0$ for all $i \in I$. Corollary \ref{PDcor} then yields $(Au^s)_r = A_{rs} - (u^r_I)^TA_Iu^s_I$, $(Au^r)_s = A_{rs} + \sum_{j \in I} u^r_j(Au^s)_j - (u^r_I)^TA_Iu^s_I$. By Lemma \ref{first_order} both these expressions are nonnegative.

It follows that
\begin{equation} \label{Delta_ineq}
(Au^s)_r + (Au^r)_s = 2(Au^s)_r + \sum_{i \in I} u^r_i(Au^s)_i \geq \sum_{i \in I} u^r_i(Au^s)_i = \sum_{i \in I} u^r_i(Au^s)_i + \sum_{i \in I} u^s_i(Au^r)_i.
\end{equation}
Switching the roles of $r,s$, we get that the inequality between the left-most and the right-most expression in \eqref{Delta_ineq} is valid also for $r \leq s$.

Let now $v$ be a zero such that $\Supp{v} \subset I \cup \{1,\dots,m\}$. Set $y = v - \sum_{r=1}^m v_ru^r$, then $\Supp{y} \subset I$ and $y_I^TA_Iy_I = v_I^TA_Iv_I - 2\sum_{r=1}^m v_rv_I^TA_Iu^r_I + \sum_{r,s=1}^m v_rv_s(u^r_I)^TA_Iu^s_I$. We obtain
\begin{eqnarray*}
0 &=& v^TAv = v_I^TA_Iv_I + 2\sum_{i \in I}\sum_{r=1}^m A_{ir}v_iv_r + \sum_{r,s = 1}^m A_{rs}v_{r}v_{s} \\ &=& v_I^TA_Iv_I + 2\sum_{i \in I}\sum_{r=1}^m \left((Au^r)_i - \sum_{j \in I}A_{ij}u^r_j\right)v_iv_r + \sum_{r,s = 1}^m \left((Au^s)_r - \sum_{j \in I} u^s_j(Au^r)_j + (u^r_I)^TA_Iu^s_I\right)v_{r}v_{s} \\ &=& y_I^TA_Iy_I + 2\sum_{i \in I}\sum_{r=1}^m (Au^r)_iv_iv_r + \frac12\sum_{r,s = 1}^m \left((Au^s)_r + (Au^r)_s - \sum_{j \in I} u^s_j(Au^r)_j - \sum_{j \in I} u^r_j(Au^s)_j\right)v_{r}v_{s} \\ &\geq& y_I^TA_Iy_I.
\end{eqnarray*}
Here for the third equality we used Corollary \ref{PDcor}, and the last inequality follows from $(Au^r)_i \geq 0$ by virtue of Lemma \ref{first_order} and \eqref{Delta_ineq}. But $A_I \succ 0$, hence $y_I = 0$ and consequently $y = 0$.

Thus $v$ is a weighted sum of the minimal zeros $u^1,\dots,u^m$ with nonnegative coefficients $\alpha_l = v_{l}$. Assume that $v$ is a minimal zero. Then by Lemma \ref{uniqueness} only one of the zeros $u^1,\dots,u^m$ can have a positive coefficient $\alpha_l$, and $v$ must be proportional to that $u^l$.
\end{proof}

Theorem \ref{PDth} restricts the ensemble of minimal zeros that a copositive matrix can have. For example, we have the following restriction on pairs of minimal zeros with overlapping supports.

{\corollary Let $A$ be a copositive matrix and $u,v$ minimal zeros of $A$ with supports $\Supp{u} = I$, $\Supp{v} = J$. Assume that $J \setminus I = \{k\}$ consists of one element. Then every zero $w$ of $A$ with support $\Supp{w} \subset I \cup J$ can be represented as a convex conic combination $w = \alpha u + \beta v$ with $\alpha,\beta \geq 0$. In particular, up to multiplication by a positive constant, there are no minimal zeros $w$ with $\Supp{w} \subset I \cup J$ other than $u$ and $v$. }

\begin{proof}
By Lemma \ref{uniqueness} there exists $i \in I$ such that $i \not\in J$. By Corollary \ref{cor:corank1} the principal submatrix $A_{I \setminus \{i\}}$ is positive definite. Noting that $J \setminus \{k\} \subset I \setminus \{i\}$, the proof is accomplished by applying Theorem \ref{PDth} to this submatrix and to the zeros $u^1 = v$, $u^2 = u$.
\end{proof}

\section{Irreducibility of copositive matrices} \label{sec:irred}

In this section we establish necessary and sufficient criteria for the irreducibility of a copositive matrix $A \in \copos{n}$ with respect to the cones $\NNM{n}$ and $S_+(n)$, respectively.

First we give a slightly stronger version of Lemma \ref{irreducibility_old}, by requiring the zero $u$ to be {\it minimal}.

{\lemma \label{irreducibility_new} Let $A \in \copos{n}$, and let $i,j \in \{1,\dots,n\}$. Then $A$ is irreducible with respect to $E_{ij}$ if and only if there exists a minimal zero $u$ of $A$ such that $(Au)_i = (Au)_j = 0$ and $u_i + u_j > 0$. }

\begin{proof}
If there exists a minimal zero with the required properties, then $A$ is irreducible with respect to $E_{ij}$ by Lemma \ref{irreducibility_old}. Hence we have to prove only the "only if" direction.

Assume that $A$ is irreducible with respect to $E_{ij}$. By Lemma \ref{irreducibility_old} there exists a zero $v$ of $A$ with $(Av)_i = (Av)_j = 0$ and $v_i + v_j > 0$. Let without restriction of generality $v_i > 0$. By Corollary \ref{zero_decomp} there exist minimal zeros $u,\dots,w$ of $A$ such that $v = u + \dots + w$. Let without loss of generality $u_i > 0$. By Lemma \ref{first_order} we have $Au \geq 0,\dots,Aw \geq 0$. From $(Av)_i = 0$ it then follows that $(Au)_i = \dots = (Aw)_i = 0$. Similarly, we obtain $(Au)_j = 0$. Thus $u$ is a minimal zero with the required properties.
\end{proof}

{\corollary \label{irreducibility_full} Let $A \in \copos{n}$. Then $A$ is irreducible with respect to $\NNM{n}$ if and only if for every pair of indices $i,j \in \{1,\dots,n\}$ there exists a minimal zero $u$ of $A$ such that $(Au)_i = (Au)_j = 0$ and $u_i + u_j > 0$. \qed }

%\medskip
%
%We can rewrite the condition given by Corollary \ref{irreducibility_full} in the following, more explicit, form.
%
%{\lemma \label{irreducibility_full0} Let $A \in \copos{n}$. Then $A$ is irreducible with respect to $\NNM{n}$ if and only if for every pair of indices $i,j \in \{1,\dots,n\}$ we have
%\[ \min \left( \{ (Au)_j \,|\, u \in \MSOZ{A},\ u_i > 0 \} \cup \{ (Au)_i \,|\, u \in \MSOZ{A},\ u_j > 0 \} \right) = 0.
%\] }
%
%\begin{proof}
%Assume that $A$ is irreducible with respect to $\NNM{n}$. By Corollary \ref{irreducibility_full}, for every index pair $i,j$ the minimum in question ranges over a nonempty set which contains at least one zero. By Lemma \ref{first_order} the minimum is nonnegative. Hence the minimum is zero.
%
%The converse direction is a direct consequence of Corollary \ref{irreducibility_full}.
%\end{proof}

\medskip

We shall now consider irreducibility with respect to the cone of positive semi-definite matrices.

{\lemma \label{irred_rank1} Let $A \in \copos{n}$ be a copositive matrix and let $w \in \mathbb R^n$ be a nonzero vector. Then $A$ is irreducible with respect to $ww^T$ if and only if there exists a zero $u$ of $A$ with $w^Tu \not= 0$. }

\begin{proof}
Let us first assume that there exists a zero $u$ with $w^Tu \not= 0$. For every $\varepsilon > 0$, we then have $u^T(A-\varepsilon ww^T)u = -\varepsilon(w^Tu)^2 < 0$, and $A - \varepsilon ww^T \not\in \copos{n}$. It follows that $A$ is irreducible with respect to $ww^T$.

It remains to show the ``only if" direction. Let $A \in \copos{n}$ be irreducible with respect to $ww^T$. For every $\varepsilon > 0$, consider the optimization problem
\begin{equation} \label{opt_prob}
\min_v \frac12 v^T(A - \varepsilon ww^T)v\ :\ v \geq 0, {\bf 1}^Tv = 1.
\end{equation}
The optimal value of this problem is negative, and it is attained by compactness of the feasible set. Let $v$ be a global minimizer of the problem. Having only linear constraints, the problem fulfills a constraint qualification \cite{Abadie67},\cite[p.52]{GSD05},
and therefore it follows from the Karush-Kuhn-Tucker optimality conditions that there exist Lagrange multipliers $\lambda \in \mathbb R_+^n$ and $\mu \in \mathbb R$ such that $v^T\lambda = 0$ and $(A - \varepsilon ww^T)v - \lambda + \mu{\bf 1} = 0$. Multiplying with $v^T$, we obtain $\mu = -v^T(A - \varepsilon ww^T)v > 0$. From $v^T\lambda = 0$ it also follows that $\Supp{\lambda} \cap \Supp{v} = \emptyset$.

Let now $\varepsilon_k \to 0$ be a sequence, let $v^k \in \mathbb R_+^n$ be a global minimizer of problem~(\ref{opt_prob}) for $\varepsilon = \varepsilon_k$, and let $\lambda^k = (\lambda_1^k,\dots,\lambda_n^k)$, $\mu^k$ be the corresponding Lagrange multipliers. Note that $\langle \lambda^k,v^k \rangle = 0$, $\mu^k > 0$, and
\begin{equation} \label{grad_general}
(A - \varepsilon_k ww^T)v^k - \lambda^k + \mu^k{\bf 1} = 0
\end{equation}
holds for all $k$.

By possibly choosing a subsequence, we can assume without restriction of generality that $v^k \to u$ for some vector $u \geq 0$ with ${\bf 1}^Tu = 1$. We may assume without loss of generality that $\Supp{u} \subset \Supp{v^k}$, and hence $\Supp{\lambda^k} \cap \Supp{u} = \emptyset$ and $u^T\lambda^k = 0$ for all $k$. Multiplying \eqref{grad_general} by $(v^k - u)^T$, we then get
\begin{equation} \label{diff_zero}
(v^k - u)^T(A - \varepsilon_k ww^T)v^k = 0.
\end{equation}

Moreover, $0 \geq \lim_{k \to \infty} (v^k)^T(A - \varepsilon_k ww^T)v^k = u^TAu \geq 0$, and $u \in \SOZ{A}$. By Lemma~\ref{first_order} we have $Au \geq 0$. Suppose the index $i$ is such that $(Au)_i > 0$. From (\ref{grad_general}) we obtain $((A - \varepsilon_k ww^T)v^k)_i = \lambda_i^k - \mu^k < \lambda_i^k$ for all $k$. By $\lim_{k \to \infty}((A - \varepsilon_k ww^T)v^k)_i = (Au)_i > 0$ we must have $\lambda_i^k > 0$ for $k$ large enough. We may hence assume without loss of generality that $\Supp{Au} \subset \Supp{\lambda^k}$ for all $k$. It follows that $\Supp{Au} \cap \Supp{v^k} = \emptyset$, and hence $(v^k)^TAu = 0$.

Inserting this into \eqref{diff_zero}, we obtain $(v^k)^T(A - \varepsilon_k ww^T)v^k = -\varepsilon_ku^Tww^Tv^k$. Hence $w^Tu \cdot w^Tv^k > 0$, and thus $w^Tu \not= 0$. Hence $u$ is the required zero.
\end{proof}

As with the case of irreducibility with respect to $E_{ij}$, we may require the zero to be minimal.

{\corollary \label{irred_min_rank1} Let $A \in \copos{n}$ be a copositive matrix and let $w \in \mathbb R^n$ be a nonzero vector. Then $A$ is irreducible with respect to $ww^T$ if and only if there exists a minimal zero $u$ of $A$ with $w^Tu \not= 0$. }

\begin{proof}
By Lemma \ref{irred_rank1}, the existence of a minimal zero $u$ with $w^Tu \not= 0$ implies irreducibility with respect to $ww^T$.

Let, on the other hand, $A$ be irreducible with respect to $ww^T$. Then by Lemma \ref{irred_rank1} there exists a zero $v$ of $A$ with $w^Tv \not= 0$. By Corollary \ref{zero_decomp} there exist minimal zeros $u,\dots,y$ such that $v = u + \dots + y$. For at least one of these minimal zeros, let it be $u$, we then must have $w^Tu \not= 0$. This concludes the proof.
\end{proof}

We are now able to characterize irreducibility with respect to the cone of positive semi-definite matrices in terms of minimal zeros.

{\theorem \label{S_irred} A copositive matrix $A \in \copos{n}$ is irreducible with respect to the cone $S_+(n)$ if and only if the linear span of the minimal zeros of $A$ equals $\mathbb R^n$. In particular, the number of linearly independent minimal zeros has to be at least $n$. }

\begin{proof}
The matrix $A$ is irreducible with respect to the cone $S_+(n)$ if and only if it is irreducible with respect to all extreme rays of this cone. By Corollary \ref{irred_rank1}, this holds if and only if for every nonzero vector $w \in \mathbb R^n$ there exists a minimal zero $u$ of $A$ such that $w^Tu \not= 0$. This condition holds if and only if the minimal zeros span the whole space.
\end{proof}

\section{Minimal support sets of irreducible copositive matrices} \label{sec:exceptional}

In this section we obtain necessary conditions for a collection ${\cal I} = \{I_1,\dots,I_m\}$ of index subsets $I_i \subset \{1,\dots,n\}$ to represent the minimal support set $\Supp\MSOZ{A}$ of a copositive matrix $A \in \copos{n}$ which is irreducible with respect to both $S_+(n)$ and ${\cal N}_n$ and satisfies $A_{ii} = 1$ for all $i$.

The obtained results can be applied to the classification of the extreme rays of $\copos{n}$. The extremal elements of $\copos{n}$ which are positive semi-definite or nonnegative have been described in \cite[Theorem 3.2]{HallNewman63}. Extremal elements which are neither positive semi-definite nor nonnegative, i.e., which are exceptional, are necessarily irreducible with respect to both $S_+(n)$ and $\NNM{n}$. Following \cite[p.9]{Baumert67} and \cite[p.1615]{DDGH13a}, we may limit our consideration to extreme elements $A \in \copos{n}$ satisfying $A_{ii} = 1$ for all $i$. The results of this section limit the number of possible minimal support sets which can occur in an exceptional extreme element with unit diagonal. The classification of the extreme rays at least for $\copos{6}$ thus comes within reach.

Before we state the main result of this section, we will need to consider the connection between the linear span of the minimal zero set $\MSOZ{A}$ and the properties of the minimal support set $\Supp{\MSOZ{A}}$. Let $I_1,\dots,I_m \subset \{1,\dots,n\}$ be the elements of ${\cal I} = \Supp{\MSOZ{A}}$, sorted by their cardinality. Since $A_{ii} = 1$ for all $i$, we cannot have a zero $u \in \SOZ{A}$ with exactly one positive element. Hence the cardinalities of the support sets satisfy $|I_k| \geq 2$, $k = 1,\dots,m$. Let $m_2$ be the number of support sets of cardinality 2.

We now construct two graphs $G_2({\cal I}),G_{>2}({\cal I})$ from $I_1,\dots,I_m$. The graph $G_2({\cal I})$ has $n$ vertices $1,\dots,n$ and $m_2$ edges $I_1,\dots,I_{m_2}$. The graph $G_{>2}({\cal I})$ is bipartite, with the two vertex subsets being defined as $V = \{1,\dots,n\}$, $W = \{m_2+1,\dots,m\}$. A pair $(v,w) \in V \times W$ is an edge of $G_{>2}({\cal I})$ if and only if $v \in I_w$. Let $G_{2,1},\dots,G_{2,r}$ be the connected components of $G_2({\cal I})$ which are bipartite.

{\lemma \label{graph_aux} Let $G$ be a graph with vertex set $\{1,\dots,n\}$ and $m$ edges $(u_1,v_1),\dots,(u_m,v_m)$. Let $M_G$ be the $n \times m$ matrix whose entries $(u_k,k)$ and $(v_k,k)$, $k = 1,\dots,m$, are equal to 1 and all other entries are equal to 0. Let $G_1,\dots,G_r$ be those connected components of $G$ which are bipartite. Then for every $k = 1,\dots,r$ the rows of $M_G$ with indices in the vertex set of $G_k$ are linearly dependent. }

\begin{proof}
Let without loss of generality $1,\dots,p$ be the vertices of $G_k$, and let $r^1 = (r^1_1,\dots,r^1_m),\dots,r^p = (r^p_1,\dots,r^p_m)$ be the corresponding rows of $M_G$. Since $G_k$ is bipartite, there exist values $\sigma_1,\dots,\sigma_p \in \{-1,+1\}$ such that for every edge $(i,j) \in G_k$ we have $\sigma_i\sigma_j = -1$. In other words, the assignment of the labels $\sigma_1,\dots,\sigma_p$ to the vertices $1,\dots,p$ realizes a 2-coloring of $G_k$.

By construction, the column $j$ of $M_G$ has exactly two entries which equal 1, and all its other entries equal to 0. The row indices of the entries which equal 1 are given by the vertices $u_j,v_j$ of edge $j$ of $G$. Since $G_k$ is a connected component of $G$, this edge either connects two vertices outside of $G_k$, or it is an edge of $G_k$. In the first case $r^i_j = 0$ for all $i = 1,\dots,p$. In the second case $r^{u_j}_j = r^{v_j}_j = 1$, $r^i_j = 0$ for $i \not= u_j,v_j$, and $\sigma_{u_j}\sigma_{v_j} = -1$. In both cases we obtain $\sum_{i=1}^p \sigma_ir^i_j = 0$.

Since this is valid for every $j = 1,\dots,m$, we finally obtain $\sum_{i=1}^p \sigma_ir^i = 0$. Thus the rows $r^1,\dots,r^p$ are linearly dependent.
\end{proof}

%Recall that the vertex partition of a connected bipartite graph is unique, because such a graph is uniquely 2-colorable \cite[p.223]{ChartrandZhang}.

{\lemma \label{lem:graph} Let $A \in \copos{n}$ be a copositive matrix with unit diagonal. Let $I_1,\dots,I_m$ be the elements of its minimal support set ${\cal I} = \Supp\MSOZ{A}$, ordered by cardinality, and let $m_2$ be the number of supports with cardinality 2. Define the two graphs $G_2({\cal I}),G_{>2}({\cal I})$ as above, and let $G_{2,1},\dots,G_{2,r}$ be the connected components of $G_2({\cal I})$ which are bipartite.

If the linear span of the minimal zero set $\MSOZ{A}$ is the whole space $\mathbb R^n$, then there exist edges $(v_1,w_1),\dots,(v_r,w_r)$ of $G_{>2}({\cal I})$ such that $v_j$ is a vertex of $G_{2,j}$ for all $j = 1,\dots,r$, and the vertices $w_1,\dots,w_r$ are mutually different. }

\begin{proof}
Assume the conditions of the lemma. Let $u^1 = (u^1_1,\dots,u^1_n)^T,\dots,u^n = (u^n_1,\dots,u^n_n)^T \in \MSOZ{A}$ be linearly independent minimal zeros of $A$, normalized such that $\max_l u^k_l = 1$ for all $k = 1,\dots,n$. For $k = 1,\dots,n$, let $i_k$ be the index of the support of $u^k$, $\Supp{u^k} = I_{i_k}$. By Lemma \ref{uniqueness} the indices $i_1,\dots,i_n$ are mutually distinct, and we may assume without loss of generality that they are ordered. Let $s$ be the number of indices $i_k$ such that $|I_{i_k}| = 2$, i.e., $|I_{i_k}| = 2$ for $k \leq s$ and $|I_{i_k}| > 2$ for $k > s$.

Then the vectors $u^1,\dots,u^s$ appear as columns of the $n \times m_2$ matrix $M_{G_2({\cal I})}$ built from the graph $G_2({\cal I})$ as in Lemma \ref{graph_aux}. The nonzero entries of the columns $u^{s+1},\dots,u^n$ correspond to edges in the graph $G_{>2}({\cal I})$.

The $n \times n$ matrix $U$ composed of the column vectors $u^1,\dots,u^n$ is invertible. Denoting by $\sign{\sigma}$ the sign of the permutation $\sigma \in S_n$, we get $\det U = \sum_{\sigma \in S_n} \sign{\sigma} \prod_{k=1}^n u^k_{\sigma(k)} \not= 0$. We shall partition the set of the summands in this sum into $\frac{n!}{s!}$ subsets indexed by the $(n-s)$-tuples $(\sigma(s+1),\dots,\sigma(n))$.

Let $J = (j_{s+1},\dots,j_n) \in \{1,\dots,n\}^{n-s}$ be an ordered $(n-s)$-tuple of mutually distinct indices, and let ${\cal J}$ be the set of all $\frac{n!}{s!}$ such $(n-s)$-tuples. For $J \in {\cal J}$, let $j_1 < \dots < j_s$ be the indices in the complement $\{1,\dots,n\} \setminus \{j_{s+1},\dots,j_n\}$. Let further $\sigma_J \in S_n$ be the permutation given by $\sigma_J(k) = j_k$, $k = 1,\dots,n$, and let $U_J$ be the $s \times s$ submatrix of $U$ with columns $1,\dots,s$ and rows $j_1,\dots,j_s$. Then the determinant can be written as $\det U = \sum_{J \in {\cal J}} \sign{\sigma_J} \det U_J \prod_{k=s+1}^n u^k_{j_k}$.

It follows that for at least one ordered $(n-s)$-tuple $J \in {\cal J}$ the summand $\det U_J \prod_{k=s+1}^n u^k_{j_k}$ is nonzero. Let without loss of generality $J = (s+1,\dots,n)$, such that $\sigma_J$ is the neutral element of $S_n$ and $U_J$ is the upper left $s \times s$ block of $U$. Otherwise we may permute the coordinates of $\mathbb R^n$, or equivalently, the rows of $U$, to achieve $J = (s+1,\dots,n)$. We then have that $U_J$ is invertible and $u^k_k \not= 0$ for $k = s+1,\dots,n$.

The columns of $U_J$ appear as columns $i_1,\dots,i_s$ in the $s \times m_2$ matrix formed of the first $s$ rows of $M_{G_2({\cal I})}$. Since the columns of $U_J$ are linearly independent, the first $s$ rows of $M_{G_2({\cal I})}$ must also be linearly independent. From Lemma \ref{graph_aux} it then follows that the vertex subset $\{1,\dots,s\}$ cannot be a superset of the vertex set of $G_{2,j}$, $j = 1,\dots,r$. Hence for every $j = 1,\dots,r$ there exists a vertex $v_j$ of $G_{2,j}$ such that $v_j > s$.

Now note that for $k > s$ we have $u^k_l \not= 0$ if and only if $(l,i_k)$ is an edge of $G_{>2}({\cal I})$. Hence $(k,i_k)$ is an edge of $G_{>2}({\cal I})$ for all $k = s+1,\dots,n$.

Then the edges $(v_j,i_{v_j})$, $j = 1,\dots,r$, witness the validity of the claim of the lemma.
\end{proof}

%Let us construct an $n \times m$ matrix ${\cal M}({\cal I})$ from $I_1,\dots,I_m$ as follows:
%\begin{equation} \label{defMA}
%({\cal M}({\cal I}))_{ij} = \left\{ \begin{array}{rcl} 0, &\quad& i \not\in I_j, \\
%1, & & i \in I_j,\ |I_j| = 2, \\
%x_{ij}, & & i \in I_j,\ |I_j| > 2, \end{array} \right.
%\end{equation}
%where $x_{ij}$ are real indeterminates. Note that for copositive matrices $A$ with unit diagonal, there cannot be supports consisting of one index only, and hence automatically $|I_j| \geq 2$ for all $j$.
%
%{\lemma \label{rank_pattern} Let $A \in \copos{n}$ be irreducible with respect to the cone $S_+(n)$, and let $A_{ii} = 1$ for all $i$. Let $I_1,\dots,I_m$ be the elements of the minimal support set ${\cal I} = \Supp\MSOZ{A}$, and let the matrix ${\cal M}({\cal I})$ be defined by \eqref{defMA}. Then there exists a realization of the indeterminates $x_i,x_{ij}$ such that ${\cal M}({\cal I})$ has rank $n$. }
%
%\begin{proof}
%By Theorem \ref{S_irred} the span of the set $\MSOZ{A}$ of minimal zeros of $A$ equals $\mathbb R^n$. For every $k = 1,\dots,m$, choose a minimal zero $u^k$ of $A$ such that $\Supp{u^k} = I_k$. If $|I_k| = 2$, then by Lemma \ref{two_zeros} we may choose $u^k$ such that its nonzero elements equal 1. By Lemma \ref{uniqueness} the linear span of $\MSOZ{A}$ equals $\spa\{u^1,\dots,u^m\}$. Hence by Theorem \ref{S_irred} the $n \times m$ matrix $U = (u^1,\dots,u^m)$ has rank $n$. By construction $U$ determines a realization of the indeterminate matrix ${\cal M}({\cal I})$.
%\end{proof}

{\lemma \label{trig_desc} Let $A \in \copos{n}$ be irreducible with respect to ${\cal N}_n$ and such that $A_{ii} = 1$ for every $i = 1,\dots,n$. Then for every $i,j = 1,\dots,n$ there exists a unique real number $\alpha_{ij} \in [0,1]$ such that $A_{ij} = -\cos(\alpha_{ij}\pi)$. }

\begin{proof}
We have $A_{ij} \geq -1$ because the principal submatrix $A_{\{i,j\}}$ has to be copositive. On the other hand, $A_{ij} \leq 1$ by \cite[Lemma 3.1]{HoffmanPereira73} because $A$ is irreducible with respect to ${\cal N}_n$. The claim of the lemma now readily follows.
\end{proof}

Clearly we have $\alpha_{ij} = \alpha_{ji}$ and $\alpha_{ii} = 1$ for all $i,j$. We shall now consider linear equalities and inequalities on these scalars imposed by the minimal support set of $A$. First we shall provide some auxiliary results.

{\lemma \label{ineq5} Let $A \in \copos{5}$ with $A_{ii} = 1$ for $i = 1,\dots,5$ be irreducible with respect to $E_{ij}$ for all $i,j = 1,\dots,5$. Let $\alpha_{ij} \in [0,1]$ such that $A_{ij} = -\cos(\alpha_{ij}\pi)$, $i,j = 1,\dots,5$. Then $\sum_{1 \leq i < j \leq 5} \alpha_{ij} \geq 4$. }

\begin{proof}
If $A \not\in S_+(5)$, then by Lemma \ref{ineq4} $A$ can be brought to the form \eqref{Smatrix} with $\sum_{i=1}^5 \theta_i < \pi$ by a permutation of its rows and columns. It follows that in this case $\sum_{1 \leq i < j \leq 5} \alpha_{ij} = 5 - \frac{1}{\pi}\sum_{i=1}^5 \theta_i > 4$.

Let now $A \in S_+(5)$, and let $B$ be defined element-wise by $B_{ij} = \frac{2}{\pi}\arcsin A_{ij} = 1 - \frac{2}{\pi}\arccos A_{ij} = 2\alpha_{ij} - 1$, $i,j = 1,\dots,5$. Then $B \in {\cal MC}_5$ by Lemma \ref{lem_maxcut}. By Definition \ref{def_maxcut} the extremal elements of the polytope ${\cal MC}_5$ are given by $vv^T$ with $v_i \in \{-1,+1\}$ for all $i = 1,\dots,5$. For every vector $v \in \{-1,+1\}^5$ we have, however, $\sum_{1 \leq i < j \leq 5} v_iv_j \geq -2$, and hence $-2 \leq \sum_{1 \leq i < j \leq 5} B_{ij} = 2\sum_{1 \leq i < j \leq 5} \alpha_{ij} - 10$. This yields the claim of the lemma.
\end{proof}

{\corollary \label{cor:ineq5} Let $A \in \copos{5}$ with unit diagonal and $-1 \leq A_{ij} \leq 1$ for all $i,j = 1,\dots,5$. Let
$\alpha_{ij} \in [0,1]$ such that $A_{ij} = -\cos(\alpha_{ij}\pi)$ for $i,j = 1,\dots,5$. Then
$\sum_{1\leq i<j\leq5} \alpha_{ij} \geq 4$. }

\begin{proof}
The matrix $A$ can be decomposed into a sum $\Gamma + N$, where $\Gamma \in \copos{5}$ is irreducible with respect to $E_{jk}$ for all $1 \leq j < k \leq 5$, and $N \in \NNM{5}$ with zero diagonal. Let $\gamma_{jk} \in [0,1]$ be such that $\Gamma_{jk} = -\cos(\gamma_{jk}\pi)$ for $1 \leq j < k \leq 5$. Then we have $\sum_{1 \leq j < k \leq 5} \gamma_{jk} \geq 4$ by Lemma \ref{ineq5}. But $\gamma_{jk} \leq \alpha_{i_ji_k}$ for all $1 \leq j < k \leq 5$, because the function $f(x) = -\cos(\pi x)$ is strictly increasing on $[0,1]$. The claim now readily follows.
\end{proof}

{\lemma \label{lin_rel} Let $A \in \copos{n}$ be irreducible with respect to ${\cal N}_n$ and such that $A_{ii} = 1$ for every $i = 1,\dots,n$. For $i,j = 1,\dots,n$, let $\alpha_{ij} \in [0,1]$ be such that $A_{ij} = -\cos(\alpha_{ij}\pi)$. Let further $B$ be a real symmetric $n \times n$ matrix defined element-wise by $B_{ij} = \frac{2}{\pi}\arcsin A_{ij} = 2\alpha_{ij} - 1$. Then the following relations hold, where the indices $i,j,k$ are assumed to be pairwise distinct:
\begin{itemize}
\item[(a)] if $\{i,j\} \in \Supp{\MSOZ{A}}$, then $\alpha_{ij} = 0$;
\item[(b)] if $\{i,j\} \not\in \Supp{\MSOZ{A}}$, then $\alpha_{ij} > 0$;
\item[(c)] if $I \in \Supp{\MSOZ{A}}$, then $B_I \in {\cal MC}_{|I|}$;
\item[(d)] if $I \subset J$ strictly and $J \in \Supp{\MSOZ{A}}$, then $B_I \in \relint{\cal MC}_{|I|}$;
%if $\{i,j\} \subset \Supp{u}$ for some $u \in \MSOZ{A}$, then $\alpha_{ij} < 1$;
\item[(e)] if $\{i,j,k\} \in \Supp{\MSOZ{A}}$, then $\alpha_{ij} + \alpha_{ik} + \alpha_{jk} = 1$;
\item[(f)] if there does not exist $I \in \Supp{\MSOZ{A}}$ such that $I \subset \{i,j,k\}$, then $\alpha_{ij} + \alpha_{ik} + \alpha_{jk} > 1$;
\item[(g)] if $\{i,j\} \in \Supp{\MSOZ{A}}$, then $\alpha_{ik} + \alpha_{jk} \geq 1$ for all $k$;
\item[(h)] for every pairwise distinct indices $i_1,\dots,i_5 \in \{1,\dots,n\}$ we have $\sum_{1 \leq j < k \leq 5}\alpha_{i_ji_k} \geq 4$.
\end{itemize} }

\begin{proof}
\begin{itemize}
\item[(a)] By Lemma \ref{corank1} we have $A_{ij} = -1$ and hence $\alpha_{ij} = 0$.
\item[(b)] If $\alpha_{ij} = 0$, then $A_{ij} = -1$ and $\{i,j\} \in \Supp{\MSOZ{A}}$ by Lemma \ref{corank1}.
\item[(c)] By Lemma \ref{corank1} the principal submatrix $A_I$ is positive semi-definite. The claim now follows from Lemma \ref{lem_maxcut}.
\item[(d)] By Lemma \ref{corank1} the principal submatrix $A_I$ is positive definite. The claim now follows from Lemma \ref{lem_maxcut} and the fact that the function $f(x) = \frac{2}{\pi}\arcsin x$ is bijective on $[-1,1]$.
%Let $\alpha_{ij} = 1$, then $A_{ij} = 1$ and the principal submatrix $A_{\{i,j\}}$ has the kernel vector $(1,-1)^T$. By Lemma \ref{corank1} this contradicts the existence of an index set $I \in \Supp{\MSOZ{A}}$ such that $\{i,j\} \subset I$\footnote{This fact has been brought to the attention of the author by Peter Dickinson, who has proven it by other means.}.
\item[(e)] By Lemma \ref{irreducibility_new} the principal submatrix $A_{\{i,j,k\}}$ is irreducible with respect to $\NNM{3}$. The relation $\alpha_{ij} + \alpha_{ik} + \alpha_{jk} = 1$ now follows from Lemma \ref{3zero}.
\item[(f)] By Lemma \ref{irreducibility_new} the principal submatrix $A_{\{i,j,k\}}$ is not irreducible with respect to $E_{12}$, $E_{13}$, and $E_{23}$. Hence it can be decomposed into a sum $\Gamma + N$, where $\Gamma \in \copos{3}$ is irreducible with respect to $\{E_{12},E_{13},E_{23}\}$, and $N \in \NNM{3}$ with zero diagonal and positive off-diagonal elements. Let $\gamma_{lm} \in [0,1]$ be such that $\Gamma_{lm} = -\cos(\gamma_{lm}\pi)$ for $1 \leq l < m \leq 3$. By Lemma \ref{3zero} we then have $\gamma_{12} + \gamma_{13} + \gamma_{23} = 1$. However, $\alpha_{ij} > \gamma_{12}$, $\alpha_{ik} > \gamma_{13}$, and $\alpha_{jk} > \gamma_{23}$, because the function $f(x) = -\cos(\pi x)$ is strictly increasing on $[0,1]$. The claim now readily follows.
\item[(g)] Let $u \in \MSOZ{A}$ with $\Supp{u} = \{i,j\}$. By Lemma \ref{two_zeros} the two positive elements of $u$ are equal, and we may normalize them to 1. By Lemma \ref{first_order} we get $(Au)_k = A_{ik} + A_{jk} = \cos((1-\alpha_{ik})\pi)-\cos(\alpha_{jk}\pi) \geq 0$. It follows that $(1-\alpha_{ik})\pi \leq \alpha_{jk}\pi$ and hence $\alpha_{ik} + \alpha_{jk} \geq 1$.
\item[(h)] This follows from Corollary \ref{cor:ineq5}.
\end{itemize}
\end{proof}

{\theorem \label{exceptional_th} Let $A \in \copos{n}$ be a copositive matrix satisfying $A_{ii} = 1$ for all $i$. Suppose that $A$ is irreducible with respect to both $S_+(n)$ and $\NNM{n}$. Let $I_1,\dots,I_m$ be the supports in the minimal support set ${\cal I} = \Supp{\MSOZ{A}}$ of $A$, ordered by their cardinality. Then ${\cal I}$ satisfies the following conditions.

\begin{itemize}
\item[(i)] Every index set $I_i$ contains $2 \leq |I_i| \leq n-2$ indices.
\item[(ii)] There do not exist $i,j$ such that $I_i \subset I_j$ strictly.
\item[(iii)] For every index set $I \subset \{1,\dots,n\}$ and indices $i,i^1,\dots,i^l,j$ satisfying the conditions
\begin{itemize}
\item $I \subset I_i$ strictly,
\item $I_{i^r} \setminus I = \{k^r\}$ consists of exactly one element for $r = 1,\dots,l$,
\item $(I_{i^r} \cap I) \subset (I_{i^{r+1}} \cap I)$ for $r = 1,\dots,l-1$,
\item $I_j \subset I \cup \{k^1,\dots,k^l\}$,
\end{itemize}
there exists $r \in \{1,\dots,l\}$ such that $j = i^r$.
\item[(iv)] Let $G_2({\cal I})$, $G_{>2}({\cal I})$ be the graphs constructed from ${\cal I}$ as in Lemma \ref{lem:graph}, and let $G_{2,1},\dots,G_{2,r}$ be the connected components of $G_2({\cal I})$ which are bipartite. Then there exist edges $(v_1,w_1),\dots,(v_r,w_r)$ of $G_{>2}({\cal I})$ such that $v_j$ is a vertex of $G_{2,j}$ for all $j = 1,\dots,r$, and the vertices $w_1,\dots,w_r$ are mutually different.
%Let the $n \times m$ matrix ${\cal M}({\cal I})$ be defined by \eqref{defMA}, then there exists a realization of the indeterminates $x_{ij}$ such that ${\cal M}({\cal I})$ has rank $n$.
\item[(v)] The system of linear equations and strict and nonstrict inequalities which is defined by (a)--(h) of Lemma \ref{lin_rel} on the variables $\alpha_{ij} = \alpha_{ji} \in [0,1]$, $1 \leq i,j \leq n$, has a solution. \end{itemize}    }

%In addition, for every pair of indices $i,j = 1,\dots,n$ there exists a minimal zero $u$ of $A$ such that $(Au)_i = (Au)_j = 0$ and $u_i + u_j > 0$. }

\begin{proof}
Since $A_{ii} \not= 0$ for all $i$, there cannot be any zero $u$ of $A$ with $|\Supp{u}| = 1$. By Lemma \ref{supp_nminus2} there cannot be a zero $u$ with $|\Supp{u}| \geq n-1$. Hence (i) holds.

Condition (ii) follows from the definition of minimality of a zero.

Condition (iii) is a consequence of Corollary \ref{cor:corank1} and Theorem \ref{PDth}.

Condition (iv) follows from Theorem \ref{S_irred} and Lemma \ref{lem:graph}.

Condition (v) follows from Lemma \ref{lin_rel}.
\end{proof}

An exceptional extremal matrix $A \in \copos{n}$ is irreducible with respect to both $S_+(n)$ and $\NNM{n}$. Hence conditions (i)--(v) of Theorem \ref{exceptional_th} are necessary conditions for the minimal support set of an exceptional extremal copositive matrix.

For given $n$ it can be checked algorithmically whether a collection $I_1,\dots,I_m \subset \{1,\dots,n\}$ of index sets satisfies conditions (i)--(v) of Theorem \ref{exceptional_th}. While this is evident for conditions (i)--(iii), we shall consider the algorithms for checking conditions (iv),(v) in more detail.

Condition (iv) can be checked by constructing a new bipartite graph $G$ from $G_{>2}({\cal I})$. To this end we manipulate the vertex subset $V = \{1,\dots,n\}$ of $G_{>2}({\cal I})$ based on the bipartite connected components $G_{2,1},\dots,G_{2,r}$ of the graph $G_2({\cal I})$. First we delete all vertices in $V$ which do not appear in one of the connected components $G_{2,1},\dots,G_{2,r}$. Then for every $j = 1,\dots,r$ we fuse the subset of $V$ which corresponds to the vertex set of $G_{2,j}$ into one vertex $j$. The so obtained graph $G$ is again bipartite with vertex subsets $\{1,\dots,r\}$ and $\{m_2+1,\dots,m\}$. Then condition (iv) is satisfied if and only if $G$ has a matching of size $r$.

Condition (v) can be checked by solving a linear program. Note that conditions (c) and (d) of Lemma \ref{lin_rel} lead to nonstrict and strict linear inequalities on $\alpha_{ij}$, respectively, because the polytope ${\cal MC}_{|I|}$ from Definition \ref{def_maxcut} and its relative interior are described by nonstrict and strict linear inequalities, respectively, on the off-diagonal elements of the matrix $B$. In order to handle the strict inequalities, we introduce a single additional slack variable $\varepsilon$. We then turn all strict inequalities into nonstrict inequalities by adding $\varepsilon$ on the smaller side. Then we maximize $\varepsilon$ with respect to the system of linear constraints provided by conditions (a)--(h) of Lemma \ref{lin_rel}. This amounts to a linear program with $\frac{n(n-1)}{2}+1$ variables $\varepsilon,\alpha_{12},\dots,\alpha_{n-1,n}$. If the program is infeasible or its optimal value is nonpositive, then condition (v) cannot hold. On the other hand, if the primal value is unbounded or the optimal value is strictly positive, then condition (v) is satisfied. Since we have to distinguish if the optimal value is zero or strictly positive, we need to solve the linear program exactly. Note that it has integer coefficients and thus a rational solution. This solution can be obtained by the simplex method.

\medskip

Two collections $I_1,\dots,I_m$ and $J_1,\dots,J_m$ satisfying conditions (i)--(v) of Theorem \ref{exceptional_th} can be considered being equivalent if there exists a permutation $\pi \in S_n$ of the indices $1,\dots,n$ such that $\{\pi(I_1),\dots,\pi(I_m)\} = \{J_1,\dots,J_m\}$. We have computed all such collections for $n \leq 7$. The number of equivalence classes is 0 for $n \leq 4$, 2 for $n = 5$, 44 for $n = 6$, and 12378 for $n = 7$. Hence $\copos{n}$ cannot have exceptional extreme rays for $n \leq 4$, which yields a quick proof of Dianandas identity $\copos{n} = S_+(n) + \NNM{n}$ for $n \leq 4$. The two equivalence classes for the case $n = 5$, with representatives $\{\{1,2\},\{2,3\},\{3,4\},\{4,5\},\{1,5\}\}$ and $\{\{1,2,3\},\{2,3,4\},\{3,4,5\},\{1,4,5\},\{1,2,5\}\}$, are realized by the Horn form \cite{HallNewman63} and the $T$-matrices from \cite{Hildebrand12a}, respectively, which indeed exhaust the types of exceptional extreme rays of $\copos{5}$. In Table \ref{Tab1}, we list one representative of each of the 44 equivalence classes for the case $n = 6$.

\begin{table}[t]
{\small
\begin{tabular}{|c|l|c|l|}
\hline
{\small No.} & {\small $\Supp{\MSOZ{A}}$} & {\small No.} & {\small $\Supp{\MSOZ{A}}$} \\
\hline
1 & \{1,2\},\{1,3\},\{1,4\},\{2,5\},\{3,6\},\{5,6\} & 23 & \{1,2,3\},\{1,2,4\},\{1,2,5\},\{1,3,6\},\{2,4,6\},\{3,4,5,6\}\!\!\! \\
2 & \{1,2\},\{1,3\},\{1,4\},\{2,5\},\{3,6\},\{4,5,6\} & 24 & \{1,2,3\},\{1,2,4\},\{1,2,5\},\{1,3,6\},\{3,4,6\},\{3,5,6\} \\
3 & \{1,2\},\{1,3\},\{1,4\},\{2,5\},\{3,5,6\},\{4,5,6\} & 25 & \{1,2,3\},\{1,2,4\},\{1,2,5\},\{1,3,6\},\{3,4,6\},\{4,5,6\} \\
4 & \{1,2\},\{1,3\},\{1,4\},\{2,5,6\},\{3,5,6\},\{4,5,6\} & 26 & \{1,2,3\},\{1,2,4\},\{1,3,5\},\{1,4,5\},\{2,3,6\},\{2,4,6\} \\
5 & \{1,2\},\{1,3\},\{2,4\},\{3,4,5\},\{1,5,6\},\{4,5,6\} & 27 & \{1,2,3\},\{1,2,4\},\{1,3,5\},\{1,4,5\},\{2,3,6\},\{3,4,6\} \\
6 & \{1,2\},\{1,3\},\{1,4,5\},\{2,4,6\},\{3,4,6\},\{4,5,6\} & 28 & \{1,2,3\},\{1,2,4\},\{1,3,5\},\{2,4,5\},\{3,4,5\},\{2,3,6\} \\
7 & \{1,2\},\{1,3\},\{2,4,5\},\{3,4,5\},\{2,4,6\},\{3,4,6\} & 29 & \{1,2,3\},\{1,2,4\},\{1,3,5\},\{2,4,5\},\{2,3,6\},\{2,5,6\} \\
8 & \{1,2\},\{1,3\},\{2,4,5\},\{3,4,5\},\{2,4,6\},\{3,5,6\} & 30 & \{1,2,3\},\{1,2,4\},\{1,3,5\},\{2,4,5\},\{3,4,6\},\{3,5,6\} \\
9 & \{1,2\},\{3,4\},\{1,3,5\},\{2,4,6\},\{1,5,6\},\{4,5,6\} & 31 & \{1,2,3\},\{1,2,4\},\{1,3,5\},\{2,4,5\},\{1,5,6\},\{2,5,6\} \\
10 & \{1,2\},\{1,3,4\},\{1,3,5\},\{2,3,6\},\{3,4,6\},\{3,5,6\} & 32 & \{1,2,3\},\{1,2,4\},\{1,3,5\},\{2,4,5\},\{1,5,6\},\{4,5,6\} \\
11 & \{1,2\},\{1,3,4\},\{1,3,5\},\{1,4,6\},\{2,5,6\},\{3,5,6\} & 33 & \{1,2,3\},\{1,2,4\},\{1,3,5\},\{2,4,5\},\{3,5,6\},\{4,5,6\} \\
12 & \{1,2\},\{1,3,4\},\{1,3,5\},\{1,4,6\},\{3,5,6\},\{4,5,6\} & 34 & \{1,2,3\},\{1,2,4\},\{1,3,5\},\{2,4,6\},\{3,5,6\},\{4,5,6\} \\
13 & \{1,2\},\{1,3,4\},\{1,3,5\},\{2,4,6\},\{3,4,6\},\{2,5,6\} & 35 & \{1,2,3,4\},\{1,2,3,5\},\{1,2,4,6\},\{1,3,5,6\},\{2,4,5,6\},\{3,4,5,6\} \\
14 & \{1,2\},\{1,3,4\},\{1,3,5\},\{2,4,6\},\{3,4,6\},\{3,5,6\} & 36 & \{1,2\},\{1,3\},\{1,4\},\{2,5\},\{4,5\},\{3,6\},\{5,6\} \\
15 & \{1,2\},\{1,3,4\},\{1,3,5\},\{2,4,6\},\{3,4,6\},\{4,5,6\} & 37 & \{1,2\},\{1,3,4\},\{1,3,5\},\{1,4,6\},\{2,5,6\},\{3,5,6\},\{4,5,6\} \\
16 & \{1,2\},\{1,3,4\},\{1,3,5\},\{2,4,6\},\{3,5,6\},\{4,5,6\} & 38 & \{1,2\},\{1,3,4\},\{1,3,5\},\{2,4,6\},\{3,4,6\},\{2,5,6\},\{3,5,6\} \\
17 & \{1,2\},\{1,3,4\},\{2,3,5\},\{3,4,5\},\{2,4,6\},\{3,4,6\} & 39 & \{1,2,3\},\{1,2,4\},\{1,2,5\},\{1,3,6\},\{1,4,6\},\{2,5,6\},\{3,5,6\} \\
18 & \{1,2,3\},\{1,2,4\},\{1,2,5\},\{1,3,6\},\{1,4,6\},\{1,5,6\} & 40 & \{1,2,3\},\{1,2,4\},\{1,2,5\},\{1,3,6\},\{1,4,6\},\{3,5,6\},\{4,5,6\} \\
19 & \{1,2,3\},\{1,2,4\},\{1,2,5\},\{1,3,6\},\{1,4,6\},\{2,5,6\} & 41 & \{1,2,3\},\{1,2,4\},\{1,2,5\},\{1,3,6\},\{2,4,6\},\{3,4,6\},\{3,5,6\} \\
20 & \{1,2,3\},\{1,2,4\},\{1,2,5\},\{1,3,6\},\{1,4,6\},\{3,5,6\} & 42 & \{1,2,3\},\{1,2,4\},\{1,2,5\},\{1,3,6\},\{2,4,6\},\{3,5,6\},\{4,5,6\} \\
21 & \{1,2,3\},\{1,2,4\},\{1,2,5\},\{1,3,6\},\{2,4,6\},\{3,4,6\} & 43 & \{1,2,3\},\{1,2,4\},\{1,2,5\},\{1,3,6\},\{1,4,6\},\{2,5,6\},\{3,5,6\},\{4,5,6\} \\
22 & \{1,2,3\},\{1,2,4\},\{1,2,5\},\{1,3,6\},\{2,4,6\},\{3,5,6\} & 44 & \{1,2,3\},\{1,2,4\},\{1,3,5\},\{1,4,5\},\{2,3,6\},\{2,4,6\},\{3,5,6\},\{4,5,6\} \\
\hline
\end{tabular}
}
\caption{Candidate minimal support sets of exceptional extreme matrices in $\copos{6}$}
\label{Tab1}
\end{table}

Finally, we give an idea of the strength of conditions (iii)--(v) in Theorem \ref{exceptional_th} by providing in Table \ref{Tab2} the number of equivalence classes of nonempty collections satisfying conditions (i) and (ii) and different combinations of conditions (iii)--(v) of Theorem \ref{exceptional_th} for different matrix sizes $n$. The lower bounds in Table \ref{Tab2} are the numbers of equivalence classes of collections where all subsets have an equal number of elements, and where there are exactly 7 subsets, respectively. The actual numbers of classes are likely one or two orders of magnitude higher, but the effort for their computation is beyond reasonable limits. From the table one sees that condition (iii) is the strongest for $n \leq 5$, but for $n \geq 6$ condition (iv) becomes dominant.

\begin{table}[t]
\begin{center}
\begin{tabular}{|l|c|c|c|c|}
\hline
combination & \multicolumn{4}{c|}{number of equivalence classes for} \\
\cline{2-5}
of conditions & $n = 4$ & $n = 5$ & $n = 6$ & $n = 7$ \\
\hline
(i),(ii) & 10 & 150 & 15933 & $> 14028724$ \\
(i),(ii),(iv),(v) & 6 & 33 & 298 & 19807 \\
(i)--(iii),(v) & 0 & 11 & 2697 & $> 157872$ \\
(i)--(iv) & 0 & 2 & 80 & 18676 \\
(i)--(v) & 0 & 2 & 44 & 12378 \\
\hline
\end{tabular}
\end{center}
\caption{Number of equivalence classes of nonempty minimal support sets satisfying different combinations of conditions}
\label{Tab2}
\end{table}

\section{Conclusions}

In this work we introduced and considered minimal zeros of copositive matrices. We established that the minimal zeros are essentially in one-to-one correspondence with the sets of indices of their positive elements (Lemma \ref{uniqueness}), which allows for a combinatorial approach to the classification of the possible combinations of minimal zeros. The main results are Theorems \ref{PDth} and \ref{S_irred}, which restrict the combinations of minimal zeros that a copositive matrix can have. The former is valid in general, while the latter holds for matrices which are irreducible with respect to the cone of positive semi-definite matrices. Lemma \ref{lin_rel} provides relations which the minimal support set imposes on the off-diagonal elements of a matrix $A \in \copos{n}$ which has unit diagonal and is irreducible with respect to both $S_+(n)$ and $\NNM{n}$. Together with trivial restrictions coming from the definition of minimality of a zero and conditions on the number of positive elements in the zeros which have been established in \cite{Baumert67,DDGH13a} they open an approach to the classification of the exceptional extreme rays of the cone $\copos{n}$ for low $n$ (Theorem \ref{exceptional_th}). Independently of the application to the classification of extreme rays, the concept of minimal zeros might prove to be a useful tool in the study of copositive matrices in general.

\section*{Acknowledgements}

The paper has benefitted enormously from the detailed comments of an anonymous referee, which we would like to gratefully acknowledge. We would also like to thank Peter J.C. Dickinson for many motivating discussions.

%\bibliography{convexity,misc,graphs}
%\bibliographystyle{plain}

\end{document}